\documentclass[12pt,reqno]{amsart}
\usepackage{graphicx, amssymb, color,pdfsync}
\usepackage[all,cmtip]{xy}
\numberwithin{equation}{section}
\usepackage{mathrsfs}
\usepackage{hyperref}

\hypersetup{
    colorlinks=true,       
    linkcolor=blue,          
    citecolor=blue,        
    filecolor=blue,      
    urlcolor=blue           
}
\usepackage{tikz}
\usetikzlibrary{matrix,arrows}
\DeclareMathOperator{\Aut}{Aut}
\usepackage[switch]{lineno}
\usepackage{yfonts}
\advance\hoffset by -0.9 truecm

\begin{document}
\newcommand{\s}{\vspace{0.2cm}}

\newtheorem{theo}{Theorem}
\newtheorem{prop}{Proposition}
\newtheorem{coro}{Corollary}
\newtheorem{lemm}{Lemma}
\newtheorem{claim}{Claim}
\newtheorem{example}{Example}
\theoremstyle{remark}
\newtheorem*{rema}{\it Remarks}
\newtheorem*{rema1}{\it Remark}
\newtheorem*{defi}{\it Definition}

\title[Nilpotent group of automorphisms of Riemann surfaces]{Nilpotent groups of automorphisms  \\of families of Riemann surfaces}
\date{}

\author{Sebasti\'an Reyes-Carocca}
\address{Departamento de Matem\'atica y Estad\'istica, Universidad de La Frontera, Avenida Francisco Salazar 01145, Temuco, Chile.}
\email{sebastian.reyes@ufrontera.cl}

\thanks{Partially supported by Fondecyt Grants 11180024, 1190991 and Redes Grant 2017-170071}

\begin{abstract} 
In this article we extend results of Zomorrodian to determine upper bounds for the order of a nilpotent group of automorphisms of a complex $d$-dimensional family of compact Riemann surfaces, where $d \geqslant 1.$ We provide conditions under which these bounds are sharp. In addition, for the one-dimensional case we construct and describe an explicit family attaining the bound for infinitely many genera. We obtain similar results for the case of $p$-groups of automorphisms.

\end{abstract}
\maketitle
\thispagestyle{empty}

\section{Introduction and statement of the results}

The classification of groups of automorphisms of compact Riemann surfaces is a classical subject of study which has attracted considerable interest ever since Hurwitz proved that the full automorphism group of a compact Riemann surface of genus $g \geqslant 2$ is finite and that its order is at most $84(g-1).$ Later, this problem acquired a new relevance when its relationship with Teichm\"{u}ller and moduli  spaces was developed.

\s

It is classically known that there are infinitely many values of $g$ for which there exists
 a compact Riemann surface of genus $g$ with automorphism group of maximal order; they are called {\it Hurwitz curves} and correspond to  branched regular covers of the projective line with three branch values, marked with 2, 3 and 7.

\s

We recall the known fact that each finite group can be realized as a group of automorphisms of a compact Riemann surface of a suitable genus. In part due to the above, an interesting problem is to study those compact Riemann surfaces whose automorphism groups share a common property and, after that, to describe among them those possessing the maximal possible number of automorphisms.

Perhaps, the most noteworthy examples concerning that are the abelian and cyclic cases. In fact, in the late nineteenth century, Wiman showed that the largest cyclic group of automorphisms of a compact Riemann surface of genus $g \geqslant 2$ has order at most $4g+2.$ Moreover, the compact Riemann surface given by the algebraic curve\begin{equation*} \label{eWiman}y^2=x^{2g+1}-1\end{equation*}shows that this upper bound is attained for each $g$. See \cite{Wi} and also \cite{Harvey1} and \cite{K1}. 

Meanwhile, as a consequence of a result due to Maclachlan, the order of an abelian group of automorphisms of a compact Riemann surface of genus $g \geqslant 2$ is at most $4g+4$; see  \cite{Mac2}. In addition, the fact that for each $g$ there exists a compact Riemann surface of genus $g$ with a group of automorphisms isomorphic to $C_2 \times C_{2g+2},$ shows that that this upper bound is attained for each value of $g.$ 

\s

Similar bounds for special classes of groups can be found in the literature in plentiful supply. For instance,   the solvable case can be found in \cite{Che} and \cite{Gro2},  the supersolvable case in  \cite{Gro3} and \cite{Z2}, the metabelian case in   \cite{ChP} and \cite{Gro4},  the metacyclic case in \cite{Sch2} and several special cases of solvable groups in \cite{Sch1}. We also refer to the survey article \cite{Gro}.

\s

By contrast, it seems that not much is known in this respect when considering complex $d$-dimensional families of compact Riemann surfaces, for $d \geqslant 1$.  Very recently, Costa and Izquierdo in \cite{CI} proved that the maximal possible order of the automorphism group of the form $ag+b,$ where $a,b \in \mathbb{Z},$ of a complex one-dimensional family of compact Riemann surfaces of genus $g \geqslant 2,$ appearing for all genera, is $4(g+1).$ In addition, they went even further by  exhibiting an explicit equisymmetric family of non-hyperelliptic compact  Riemann surfaces attaining this bound for all $g$ (c.f. \cite{BCI}). Later, the analogous problem for complex low dimensional families ($d \leqslant 4$) was addressed in \cite{IRCR} and \cite{RR1}.


\s

The aim of this article is to deal with nilpotent groups and $p$-groups of automorphisms of complex $d$-dimensional families of compact Riemann surfaces, where $d \geqslant 1.$

We recall that the Jacobian variety $JC$ of a compact Riemann surface $C$ of genus $g$ is an irreducible principally polarized abelian variety of dimension $g.$ The relevance of the Jacobian variety lies, in part, in the classical Torelli's theorem, which establishes that $$C_1 \cong C_2 \,\, \mbox{ if and only if }\,\, JC_1 \cong JC_2.$$

In this paper we shall also consider isogenous decompositions of Jacobian varieties of certain compact Riemann surfaces with a nilpotent group of automorphisms. 

\s

{\bf Nilpotent groups acting on families of Riemann surfaces} 

\s

In \cite{MMBB} Macbeath considered homomorphisms from co-compact Fuchsian groups onto finite nilpotent groups. Since every finite nilpotent group is isomorphic to the direct product of its Sylow subgroups, after introducing the concept of {\it $p$-localization} of groups, he succeeded in providing necessary and sufficient conditions under which a given signature appears as the signature of the action of a nilpotent group of automorphisms on a compact Riemann surface. 

Soon after and based on the aforementioned Macbeath's result, Zomorrodian in \cite{Z1} proved that  the order of a nilpotent group of automorphisms of a compact Riemann surface of genus $g \geqslant 2$ is at most $16(g-1).$ Moreover, he noticed that if the previous bound is sharp then $g-1$ is a power of two and the signature of the action is $(0; 2,4,8).$

\s

Here we extend the previous result from (zero-dimensional families of) compact Riemann surfaces to  $d$-dimensional families of compact Riemann surfaces, where $d \geqslant 1.$

\begin{theo} \label{nil1} Let $d \geqslant 1, g\geqslant 2$  be integers.  Let $G$ be a nilpotent group of automorphisms of a complex $d$-dimensional family of compact Riemann surfaces $C$ of genus $g$. 
\begin{enumerate}
\item The order $G$  is at most \begin{displaymath} M_{2,d}= \left\{ \begin{array}{ll}
\,\,\,\,8(g-1) & \textrm{if $d=1$}\\
\tfrac{4}{d-1}(g-1) & \textrm{if $d \geqslant 2.$}
  \end{array} \right.
\end{displaymath} 

\item  The order of $G$ is $M_{2,d}$  if and only if the signature of the action of $G$ on $C$  is\begin{displaymath}\sigma_{2,d}=  \left\{ \begin{array}{ll}
(0; 2,2,2,4) & \textrm{if $d=1$}\\
(0; 2, \stackrel{d+3}{\ldots},2) & \textrm{if $d \geqslant 2.$}
  \end{array} \right.
\end{displaymath} 
\item If the order of $G$ is $M_{2,d}$ then $G$ is a $2$-group. In particular if, in addition, $d=1$ or $d-1$ a power of two then $g-1$ is a power of two.
\end{enumerate}
\end{theo}

\s

If $g-1$ is a power of two, in \cite{Z1} it was also proved the existence of at least one compact Riemann surface of genus $g$ with a nilpotent group of automorphisms of order $16(g-1)$, showing that this upper  bound is attained for infinitely many values of $g.$

\s

Note that for $d=2$ the previous theorem guarantees that, if the order of $G$ is maximal then $g-1$ is a power of two. We notice that the converse is also true. Indeed, following \cite{RR1},  for each $g \geqslant 2,$ there exist a complex two-dimensional family of compact Riemann surfaces of genus $g$ with a dihedral group of automorphisms of order  $$M_{2,2}=4(g-1) \,\, \mbox{ acting with signature } \,\, \sigma_{2,2}=(0; 2,2,2,2,2).$$ Thus, in particular, if $g-1$ is a power of two then the involved dihedral group is nilpotent and therefore the upper bound $M_{2,2}$ is attained. 

\s

It is worth pointing out  here that Zomorrodian's method to prove the existence of a compact Riemann surface of genus $g$ with a nilpotent group of automorphisms of order $16(g-1)$ is based on an inductive argument and does not provide neither the Riemann surface nor the nilpotent group in an explicit manner; see \cite[p. 254]{Z1}.  

\s

The following theorem shows that the upper bound $M_{d,1}$  is sharp for infinitely many values of $g$. In contrast with the zero-dimensional case, our strategy is to construct a complex one-dimensional family in an explicit enough way in order to provide a  detailed description of the family. We include an  isogeny decomposition of the associated family of Jacobian varieties.

\begin{theo} \label{nil2} For each integer  $n \geqslant 5$ there is a complex one-dimensional closed family  of compact Riemann surfaces $C$ of genus $1+2^{n-3}$ with a nilpotent group of automorphisms $G$ of order $2^n$ isomorphic to the semidirect product $$(C_2 \times \mathbf{D}_{2^{n-3}}) \rtimes C_2$$presented in terms of generators $a,b,r,s$ and relations $$r^{2^{n-3}}=s^2=(sr)^2=a^2=b^2=1, [s,b]=[r,b]=1, ara=r^{-1}, asa=sr, aba=br^{2^{n-4}}$$acting on $C$ with signature $(0; 2,2,2,4).$ Furthermore:
\begin{enumerate}
\item the family consists of at most $2^{2n-6}$ equisymmetric strata,
\item up to possibly finitely many exceptions, $C$ is non-hyperelliptic and its  automorphism group agrees with $G$, and
\item the Jacobian variety $JC$ of $C$ contains an elliptic curve isogenous to $JC_{\langle r \rangle}$ and decomposes, up to isogeny, as $$JC \sim JC_{\langle s \rangle} \times JC_{\langle b \rangle},$$where the dimensions of $JC_{\langle s \rangle}$ and $JC_{\langle b \rangle}$ are $2^{n-4}$ and $2^{n-4}+1$ respectively.
\end{enumerate}
\end{theo}

\begin{rema} \mbox{}
\begin{enumerate}
\item The cases $n=3$ and $n=4$ are exceptional in the sense that the upper bound is attained by a group with a different algebraic structure. Concretely
\begin{enumerate} 
\item for $n=3$ ($g=2$) the  bound is attained by $\mathbf{D}_4,$ and
\item for $n=4$ ($g=3$) the bound is attained by $C_2 \times \mathbf{D}_4$ and by $(C_2 \times C_4) \rtimes C_2$.\end{enumerate}See \cite{bci} and \cite{Brou}.
\item We announce that for each odd integer $d \geqslant 3,$ the  bound $M_{2,d}$ is attained for infinitely many genera. We shall deal with this problem in a forthcoming paper.
\end{enumerate}

\end{rema}

\s

{\bf $p$-groups acting on families of Riemann surfaces.} 

\s

The fact that nilpotent groups of automorphisms of  compact Riemann surfaces of maximal order turn out to be 2-groups led Zomorrodian to ask for similar bounds for the class of $p$-groups. Indeed, he proved in \cite{Z3} that if $G$ is a $p$-group of automorphisms of a compact Riemann surface of genus $g \geqslant 2$ then \begin{equation} \label{vi}|G| \leqslant \epsilon(g-1) \,\, \mbox{ where } \,\, \epsilon=\left\{ \begin{array}{ll}
\,16 & \textrm{if $p= 2$}\\
\,\,\,9 & \textrm{if $p= 3$}\\
\tfrac{2p}{p-3} & \textrm{if $p \geqslant 5,$}\end{array} \right.
\end{equation}and that  \eqref{vi} turns into an equality  if and only if the signature of the action is $$(0; 2,4,8), \,\,  (0; 3,3,9)\,\, \mbox{ and } \,\,(0; p,p,p)$$respectively.  Furthermore, in the same paper it was also proved the existence of a $p$-group of order $p^n$ acting on a compact  Riemann surface of genus $1+p^{n}/\epsilon$ for each $n \geqslant 4$, showing that the bounds \eqref{vi} are sharp for infinitely many values of $g.$

\s

The following result is a direct consequence of Theorems \ref{nil1} and \ref{nil2}.
\begin{coro}
Let $d \geqslant 1$ and $g\geqslant 2$  be integers.  If $G$ is $2$-group of automorphisms of a complex $d$-dimensional family of compact Riemann surfaces $C$ of genus $g$ then:
\begin{enumerate}
\item the order $G$ is at most  $M_{2,d}$,
\item the order of $G$ is $M_{2,d}$  if and only if the signature of the action is $\sigma_{2,d}$, and 
\item the upper bound $M_{2,1}$ is attained for infinitely many values of $g.$
\end{enumerate}
\end{coro}

The following theorem extends both the previous corollary from $p=2$ to odd prime numbers $p \geqslant 3$ and the results in \cite{Z3}  from the zero-dimensional situation  to complex $d$-dimensional families. For each rational number $t \geqslant 0$ we denote its integer part by $[t].$ 

\begin{theo}\label{ppp} Let $d \geqslant 1$ and $g \geqslant 2$ be integers and let $p \geqslant 3$ be a prime number. Let $G$ be a $p$-group of automorphisms of a complex $d$-dimensional family of compact Riemann surfaces $C$ of genus $g$. 
\s

Assume $p=3.$
\s
\begin{enumerate}
\item[(1)] The order of $G$ is at most $$M_{3,d}= \tfrac{3}{d}(g-1).$$
\item[(2)] The order of $G$ is $M_{3,d}$ if and only if the signature of the action of $G$ on $C$ is 
$$\sigma_{3,d,h}=(h; 3, \stackrel{d+3-3h}{\ldots},3) \,\, \mbox{ for some } \,\, h \in \{0, \ldots, [\tfrac{d}{3}+1]\}.$$
\end{enumerate}

\s

Assume $p \geqslant 5.$

\s

\begin{enumerate}
\item[(3)] Let $\lambda_d$ be the smallest non-negative representative of $d$ modulo 3. 
The order of $G$ is at most $$ M_{p,d}=\tfrac{2}{N}(g-1) \,\, \mbox{ where }\,\, N=\tfrac{2}{3}d + \lambda_d(\tfrac{1}{3}-\tfrac{1}{p}).$$

\item[(4)] The order of $G$ is $M_{p,d}$ if and only if the signature of the action of $G$ on $C$ is 
$$\sigma_{p,d}=(\hat{h}; p, \stackrel{d+3-3\hat{h}}{\ldots},p) \,\, \mbox{ where }\,\, \hat{h}=[\tfrac{d}{3}+1].$$
\end{enumerate}
\end{theo}

The previous theorem applied to $d=1$ says that if $p \geqslant 3$ is a prime number and if $G$ is a $p$-group of automorphism of a complex one-dimensional family of compact Riemann surfaces of genus $g$ then \begin{equation}\label{vi4}|G| \leqslant M_{p,1}= \tfrac{2p}{p-1}(g-1)\end{equation}and  the that equality holds  if and only if the signature of the action is $(1; p)$ for $p \geqslant 5,$ and $(1;3)$ or $(0; 3,3,3,3)$ for $p=3.$

\s

The following theorem provides a detailed description of a complex one-dimensional family of compact Riemann surfaces whose existence shows that the bound \eqref{vi4} is sharp for each prime $p \geqslant 3$ and for infinitely many values of $g.$

\begin{theo} \label{pexistencia} Let $p \geqslant 3$ be a prime number. For each integer $n \geqslant 3$ there is a complex one-dimensional closed family  of compact Riemann surfaces $C$ of genus $$1+\tfrac{(p-1)p^{n-1}}{2}$$ with a $p$-group of automorphisms $G$ of order $p^n$ isomorphic to the semidirect product $$C_{p^{n-1}} \rtimes_p C_{p}=\langle a,b : a^{p^{n-1}}=b^{p}=1, bab^{-1}=a^{r} \rangle,$$where $r=p^{n-2}+1,$ acting on $C$ with signature $(1; p).$ In addition, 
\begin{enumerate}
\item the family consists of  $p-1$ equisymmetric strata,
\item $C$ is elliptic-$p$-gonal, 
\item up to possibly finitely many exceptions, the automorphism group of $C$ agrees with $G$, and 
\item the Jacobian variety $JC$ of $C$ decomposes, up to isogeny, as $$JC \sim E \times A^p,$$where $E$ is an elliptic curve isogenous to $JC_G$ and $A$ is an abelian subvariety of $JC$ of dimension  $\tfrac{(p-1)p^{n-2}}{2} .$
\end{enumerate}
\end{theo}

\begin{rema1}
The groups involved in this paper have order of the form $\rho(g-1)$ where $\rho \in \mathbb{Q}$ and $g-1$ is a power of a prime number. We remark that this situation differs radically from the case in which $\rho \in \mathbb{Z}$ and $g-1$ is prime; see \cite{BJ}, \cite{IJRC}, \cite{IRC} and \cite{RR1}.
\end{rema1}

This paper is organized as follows. Section \S\ref{s2} will be devoted to briefly review the basic background:  Fuchsian groups, group actions on Riemann surfaces, the equisymmetric stratification of the moduli space and the decomposition of Jacobian varieties  with group action.  The proofs of the theorems will be given in Sections \S\ref{s3}, \S\ref{s4}, \S\ref{s5} and \S\ref{s6}.

\section{Preliminaries}\label{s2}

\subsection{Fuchsian groups} A {\it Fuchsian group} is a discrete group of automorphisms of $$\mathbb{H}=\{z \in \mathbb{C}: \mbox{Im}(z) >0 \}.$$  

If $\Delta$ is a Fuchsian group and the orbit space $\mathbb{H}_{\Delta}$ given by the action of $\Delta$ on $\mathbb{H}$ is  compact, then the algebraic structure of $\Delta$ is determined by its {\it signature}: \begin{equation} \label{sig} \sigma(\Delta)=(h; m_1, \ldots, m_l),\end{equation}where $h$ is the genus of the quotient $\mathbb{H}_{\Delta}$ and $m_1, \ldots, m_l$ are the branch indices in the universal canonical projection $\mathbb{H} \to \mathbb{H}_{\Delta}.$ The signature \eqref{sig} is called {\it degenerate} if $$h=0 \, \mbox{ and } \, l=1 \,\,\,\, \mbox{ or } \,\,\,\, h=0 \, \mbox{ and } l=2 \, \mbox{ with } \, m_1 \neq m_2.$$

Let $\Delta$ be a Fuchsian group of signature \eqref{sig}. Then\begin{enumerate}
\item $\Delta$ has a canonical presentation with generators $\alpha_1, \ldots, \alpha_{h}$, $\beta_1, \ldots, \beta_{h},$ $ \gamma_1, \ldots , \gamma_l$ and relations
\begin{equation}\label{prese}\gamma_1^{m_1}=\cdots =\gamma_l^{m_l}=\Pi_{i=1}^{h}[\alpha_i, \beta_i] \Pi_{i=1}^l \gamma_i=1,\end{equation}where $[u,v]$ stands for the commutator $uvu^{-1}v^{-1}.$
\item The elements of $\Delta$ of finite order are conjugate to powers of $\gamma_1, \ldots, \gamma_l.$
\item The Teichm\"{u}ller space of $\Delta$ is a complex analytic manifold homeomorphic to the complex ball of dimension $3h-3+l$.
\item The hyperbolic area of each fundamental region of $\Delta$ is given by $$\mu(\Delta)=2 \pi [2h-2 + \Sigma_{i=1}^l(1-\tfrac{1}{m_i})].$$ 
\item The Euler characteristic of the signature $\sigma(\Delta)$ is the rational number $$\chi(\sigma(\Delta))=-\tfrac{1}{2\pi}\mu(\Delta).$$
\end{enumerate}
We refer to the classical articles \cite{Harvey} and \cite{singerman} for further details.

Let $\Gamma$ be a group of automorphisms of $\mathbb{H}.$ If $\Delta$ is a subgroup of $\Gamma$ of finite index then $\Gamma$ is also Fuchsian and their hyperbolic areas are related by the Riemann-Hurwitz formula $$\mu(\Delta)= [\Gamma : \Delta] \cdot \mu(\Gamma).$$

\subsection{Group actions on Riemann surfaces and localization} Let $C$ be a compact Riemann surface of genus $g \geqslant 2$ and let $\mbox{Aut}(C)$ denote its  automorphism group. A finite group $G$ acts on $C$ if there is a group monomorphism $G\to \Aut(C).$ The space of orbits $C_G$ of the action of $G$ on $C$ is naturally endowed with a Riemann surface structure such that the canonical projection $C \to C_G$ is holomorphic. 

\s

By the classical uniformization theorem, there is a unique, up to conjugation, Fuchsian group $\Gamma$ of signature $(g; -)$ such that $C \cong \mathbb{H}_{\Gamma}.$ Moreover, $G$ acts on $C$ if and only if there is a Fuchsian group $\Delta$ containing $\Gamma$ together with a group  epimorphism \begin{equation*}\label{epi}\theta: \Delta \to G \, \, \mbox{ such that }  \, \, \mbox{ker}(\theta)=\Gamma.\end{equation*}

In such a case, the group $G$ is said to act on $C$ with signature $\sigma(\Delta)$ and the action is said to be {\it  represented by the surface-kernel epimorphism} $\theta.$ See \cite{Harvey}, \cite{yoibero} and \cite{singerman}

If $G$ is a subgroup of $G'$ then the action of $G$ on $C$ is said to {\it extend} to an action of $G'$ on $C$ if:\begin{enumerate}
\item there is a Fuchsian group $\Delta'$ containing $\Delta,$ 
\item the Teichm\"{u}ller spaces of $\Delta$ and $\Delta'$ have the same dimension, and
\item there exists an epimorphism $$\Theta: \Delta' \to G' \, \, \mbox{ in such a way that }  \, \, \Theta|_{\Delta}=\theta  \mbox{ and } \ker(\theta)=\ker(\Theta).$$
\end{enumerate} 

An action is called {\it maximal} if it cannot be extended in the afore introduced sense. A complete list of signatures of pairs of Fuchsian groups $\Delta$ and $\Delta'$ for which it may be possible to have an extension as before was provided by Singerman in \cite{singerman2}. 

\s

 Let $\Delta$ be a Fuchsian group of signature \eqref{sig} and let $p$ be a prime number. Define $e_i$ as the largest integer such that $p^{e_i}$ is a divisor of $m_i.$ Following  \cite{MMBB}, the signature $$\sigma_p:=(h; p^{e_1}, \ldots, p^{e_l}),$$ where the $(i+1)$-entry is dropped if $e_i=0,$ is called the {\it $p$-localization} of $\sigma=\sigma(\Delta).$ The signature $\sigma$ is called {\it nilpotent-admissible} if $\sigma_p$ is non-degenerate for each prime $p.$

\s

Macbeath proved that if $\sigma$ is a nilpotent-admissible signature then there exists at least one nilpotent group acting as a group of automorphisms of a compact Riemann surface with signature $\sigma.$ Furthermore, if in addition the signature satisfies that  $\chi(\sigma_p) \leqslant 0$ for al least one prime $p,$ then there are infinitely many nilpotent groups with the same property. See \cite[Theorem (8.1)]{MMBB} and  \cite[Theorem (8.2)]{MMBB}.

\subsection{Equisymmetric stratification} \label{strati} Let $\text{Hom}^+(C)$ denote the group of orientation preserving self-homeomorphisms of $C.$ Two actions $\psi_i: G \to \mbox{Aut}(C)$ of $G$ on $C$ are  {\it topologically equivalent} if there exist $\omega \in \Aut(G)$ and $f \in \text{Hom}^+(C)$ such that
\begin{equation}\label{equivalentactions}
\psi_2(g) = f \psi_1(\omega(g)) f^{-1} \hspace{0.5 cm} \mbox{for all} \,\, g\in G.
\end{equation}

Each homeomorphism $f$ satisfying \eqref{equivalentactions} yields an automorphism $f^*$ of $\Delta$ where $\mathbb{H}_{\Delta} \cong C_G$. If $\mathscr{B}$ is the subgroup of $\mbox{Aut}(\Delta)$ consisting of them, then $\mbox{Aut}(G) \times \mathscr{B}$ acts on the set of epimorphisms defining actions of $G$ on $C$ with signature $\sigma(\Delta)$ by $$((\omega, f^*), \theta) \mapsto \omega \circ \theta \circ (f^*)^{-1}.$$  

Two epimorphisms $\theta_1, \theta_2 : \Delta \to G$ define topologically equivalent actions if and only if they belong to the same $(\mbox{Aut}(G) \times \mathscr{B})$-orbit (see \cite{bci}, \cite{Brou},  \cite{Harvey} and \cite{McB2}). 

We remark that if the genus of $C_G$ is one then 
$\mathscr{B}$ contains  the transformations $$A_{1,n}: \alpha_1 \mapsto \alpha_1, \,\, \beta_1 \mapsto \beta_1 \alpha_1^n, \,\,  \gamma_j \to \gamma_j, \,\, \mbox{ and }\,\, A_{2,n}: \alpha_1 \mapsto \alpha_1 \beta_1^n, \,\, \beta_1 \mapsto \beta_1, \,\,  \gamma_j \to \gamma_j$$for each  $n \in \mathbb{Z}.$ See \cite[Proposition 2.5]{Brou}.


\s

Let $\mathscr{M}_g$ denote the moduli space of compact Riemann surfaces of genus $g \geqslant2.$ It is well-known  that $\mathscr{M}_g$ is endowed with a structure of complex analytic space of dimension $3g-3,$ and that for $g \geqslant4$ its singular locus $\mbox{Sing}(\mathscr{M}_g)$ agrees with the set of points representing compact Riemann surfaces with non-trivial automorphisms. 

\s

Following \cite{b}, the singular locus of $\mathscr{M}_g$ admits an {\it equisymmetric stratification} $$\mbox{Sing}(\mathscr{M}_g)= \cup_{G, \theta} \bar{\mathscr{M}}_g^{G, \theta}$$ where 
each {\it equisymmetric stratum} $\mathscr{M}_g^{G, \theta}$, if nonempty, corresponds to one topological class of maximal actions (see also \cite{Harvey}). More precisely:

\begin{enumerate}

\item the {\it equisymmetric stratum} ${\mathscr{M}}_g^{G, \theta}$ consists of  those Riemann surfaces $C$ of genus $g$ with (full) automorphism group isomorphic to $G$ such that the action is topologically equivalent to $\theta$,

\item the closure $\bar{\mathscr{M}}_g^{G, \theta}$ of  ${\mathscr{M}}_g^{G, \theta}$ is a closed irreducible algebraic subvariety of $\mathscr{M}_g$ and consists of those Riemann surfaces $C$ of genus $g$ with a group of automorphisms isomorphic to $G$ such that the action is  topologically equivalent to $\theta$, and

\item  if the equisymmetric stratum ${\mathscr{M}}_g^{G, \theta}$ is nonempty then it is a smooth, connected,
locally closed algebraic subvariety of $\mathscr{M}_{g}$ which is Zariski dense in
$\bar{\mathscr{M}}_g^{G, \theta}.$ 
\end{enumerate}

\s

In this article we employ use the following terminology.

\begin{defi} The subset of $\mathscr{M}_g$ consisting of those compact Riemann surfaces $C$ of genus $g$ with action of a given group $G$ with a given signature will be called a {\it (closed)  family}.  \end{defi}

The complex dimension of the family is the complex dimension of the Teichm\"{u}ller space associated to a Fuchsian group $\Delta$ such that $C_G \cong \mathbb{H}_{\Delta}.$ Note that  the interior of a family consists of those Riemann surfaces whose full automorphism group is isomorphic to $G$ and  is formed by finitely many equisymmetric strata which are in correspondence with the pairwise non-equivalent topological actions of $G.$  Besides, the members of the family that do not belong to the interior is formed by those Riemann surfaces that have strictly more automorphisms than $G.$

\s

\subsection{Decomposition of Jacobians with group action} \label{jacos} It is  well-known that if $G$ acts on a compact Riemann surface $C$ then this action  induces a $\mathbb{Q}$-algebra homomorphism $$\Phi : \mathbb{Q} [G] \to \mbox{End}_{\mathbb{Q}}(JC)=\mbox{End}(JC) \otimes_{\mathbb{Z}} \mathbb{Q},$$from the rational group algebra of $G$ to the rational endomorphism algebra of $JC.$

 For each $ \alpha \in {\mathbb Q}[G]$ we define the abelian subvariety $$A_{\alpha} := {\textup Im} (\alpha)=\Phi (n\alpha)(JC) \subset JC$$where $n$ is some positive integer chosen in such a way that $n\alpha \in {\mathbb Z}[G]$.

 Let  $W_1, \ldots, W_r$ be the rational irreducible representations of $G$. For each $W_j$ we denote by $V_j$ a complex irreducible representation of $G$ associated to it.  The decomposition of $1$ as the sum $e_1 + \cdots + e_r,$ where $e_j  \in \mathbb{Q}[G]$ is a uniquely determined central idempotent computed explicitly from $W_j$, yields an isogeny $$JC \sim A_{e_1} \times \cdots \times A_{e_r}$$
which is $G$-equivariant; see \cite{l-r}. Additionally, there are idempotents $f_{j1},\dots, f_{jn_j}$ such that $e_j=f_{j1}+\dots +f_{jn_j}$  where  $n_j=d_{V_j}/s_{V_j}$ is the quotient of the degree $d_{V_j}$ of $V_j$ and its Schur index $s_{V_j}$.  These idempotents provide $n_j$ subvarieties of $JC$ which are pairwise isogenous; let $B_j$ be one of them, for every $j.$ Thus, we obtain the following isogeny
\begin{equation} \label{eq:gadec}
JC \sim_G B_{1}^{n_1} \times \cdots \times B_{r}^{n_r} 
\end{equation}
called the {\it group algebra decomposition} of $JC$ with respect to $G$. See \cite{cr} and also \cite{RCR}.

If $W_1(=V_1)$ denotes the trivial representation of $G$ then $n_1=1$ and $B_{1} \sim JC_G$.

\s

Let $H$ be a subgroup of $G$ and consider the associated regular covering map $\pi_H:C \to C_H.$ It was proved in \cite{cr} that \eqref{eq:gadec} induces the isogeny  \begin{equation} \label{indjaco}JC_H \sim  B_{1}^{{n}_1^H} \times \cdots \times B_{r}^{n_r^H} \,\,\, \mbox{ with } \,\,\, {n}_j^H=d_{V_j}^H/s_{V_j}
\end{equation}where $d_{V_j}^H$ is the dimension of the vector subspace $V_j^H$ of $V_j$ of elements fixed under $H.$ 

\s

Assume that \eqref{sig} is the signature of the action of $G$ on $C$  and that this action is represented by $\theta: \Delta \to G,$ with $\Delta$ as in \eqref{prese}. Following \cite[Theorem 5.12]{yoibero} 
\begin{equation}\label{uuaa}
\dim (B_{j})=k_{V_j}[d_{V_j}(\gamma -1)+\frac{1}{2}\Sigma_{k=1}^l (d_{V_j}-d_{V_j}^{\langle \theta(\gamma_k) \rangle} )]  \,\, \mbox{ for }\,\, 2 \leqslant j \leqslant r\end{equation} where $k_{V_j}$ is the degree of the extension $\mathbb{Q} \leqslant L_{V_j}$ with $L_{V_j}$ denoting a minimal field of definition for $V_j.$

\s

The decomposition of Jacobian varieties with group actions has been extensively studied, going back to contributions of Wirtinger, Schottky and Jung. For decompositions of Jacobians with respect to special groups, we refer to \cite{CLR},   \cite{d1}, \cite{CRC},  \cite{Do}, \cite{nos}, \cite{IJR}, \cite{IRC},  \cite{KR}, \cite{PA}, \cite{d3}, \cite{yo-racsam}, \cite{yojpaa}, \cite{kanirubiyo} and \cite{Ri}.
 
\subsection*{\it Notation} We denote the cyclic group of order $n$ by $C_n$ and the dihedral group of order $2n$ by $\mathbf{D}_n.$

\section{Proof of Theorem \ref{nil1}} \label{s3}
Let $d \geqslant 1$ and $g \geqslant 2$ be integers.   We assume that $G$ is a nilpotent group acting as a group of automorphisms of a complex $d$-dimensional family  of compact Riemann surfaces $C$ of genus $g,$   and that the signature of the action of $G$ on $C$ is $\sigma=(h; m_1, \ldots, m_l).$ 

\s

Assume that $d \geqslant 2.$ Note that, as each $m_i \geqslant 2,$ the hyperbolic area $\mu$ of a fundamental domain of a Fuchsian group of signature $\sigma$ satisfies $$\mu = 2 \pi (2h-2+\Sigma_{i=1}^l\tfrac{2}{m_i}) \geqslant 2 \pi(\tfrac{h}{2}+\tfrac{d-1}{2}) \geqslant 2 \pi \tfrac{d-1}{2}.$$ Thus,  by the Riemann-Hurwitz formula, one easily obtains that$$2(g-1) = \tfrac{\mu}{2 \pi} |G| \geqslant  \tfrac{d-1}{2}|G| \,\, \iff \,\, |G| \leqslant \tfrac{4}{d-1}(g-1)$$as claimed. Now, if we assume that  $$|G| = \tfrac{4}{d-1}(g-1)\,\, \mbox{ then }\,\, \Sigma_{i=1}^l \tfrac{1}{m_i}=\tfrac{d+3}{2}-h,$$which is at most $\tfrac{l}{2}.$ It follows that $$h=0  \,\, \mbox{ and } \,\, \Sigma_{i=1}^{d+3}\tfrac{1}{m_i}=\tfrac{d+3}{2}.$$The unique solution of the equation above is $m_i=2$ for each $i,$ and then $\sigma=(0; 2, \stackrel{d+3}{\ldots},2).$ 

\s

Assume that $d = 1.$ We have only two cases to consider; namely, $(h,l)=(1,1)$ and $(h,l)=(0,4).$ In the former case it is clear that  
$\mu \geqslant \pi.$ Assume $\sigma=(0; m_1, m_2, m_3, m_4)$  and denote by $v$ the number of periods $m_i$ that are equal to 2. Note that $v \leqslant 3$ because if $v=4$ then $\mu=0.$

\begin{enumerate}
\item[(a)] If $v=0$ then each $m_i \geqslant 3$ and therefore  $\mu \geqslant \tfrac{4\pi}{3}.$ 

\s

\item[(b)] If $v=1$ then $\sigma=(0; 2,m_2,m_3,m_4)$ where $m_i \geqslant 3.$ Note that if $m_2, m_3, m_4$ were equal to 3 then the $2$-localization of $\sigma$ would be degenerate. Then, we can assume $m_4 \geqslant 4$ and therefore $\mu  \geqslant \tfrac{7\pi}{6}.$ 

\s

\item[(c)] If $v=2$ then $\sigma=(0; 2,2,m_3,m_4)$ where $m_3, m_4 \geqslant 3$ and $\mu \geqslant \tfrac{2\pi}{3}.$

\s

\item[(d)] If $v=3$ then $\sigma=(0; 2,2,2,m_4)$ where $m_4 \geqslant 3.$ Note that $m_4$ must be a power of two, since otherwise the $p$-localization of $\sigma$ would be degenerate for some prime $p \geqslant 3.$ Thus $\mu \geqslant \tfrac{\pi}{2}.$
\end{enumerate}

\s

All the above ensures that $\mu \geqslant \tfrac{\pi}{2}$ and therefore, by the Riemann-Hurwitz formula,$$2(g-1) = |G| \tfrac{\mu}{2 \pi} \geqslant  \tfrac{|G|}{4} \,\, \iff \,\, |G| \leqslant 8(g-1)$$as claimed. Now, if $|G|=8(g-1)$ then$$\Sigma_{i=1}^{l}\tfrac{1}{m_i}=\tfrac{7}{4}-h \leqslant \tfrac{4-3h}{2} \,\, \mbox{ and therefore } h=0 \,\, \mbox{ and } \,\, \Sigma_{j=1}^4 \tfrac{1}{m_j}=\tfrac{7}{4},$$showing that $m_1=m_2=m_3=2$ and $m_4=4.$ Thus, $\sigma=(0; 2,2,2,4)$ as desired. 

\s

Finally, as the group $G$ is assumed to be nilpotent and as, in each case, the genus of the corresponding quotient is zero, we can apply \cite[Theorem 2.11]{Z1} to ensure that the prime factors of $|G|$ are necessarily contained in the set of prime factors of the periods of $\sigma.$ Thus, the group $G$ is a $2$-group. Consequently, if we assume that, in addition, $d=1$ or $d-1$ is a power of two, then  we can conclude that  $g-1$ is a power of two as well.

\section{Proof of Theorem 2} \label{s4}

Let $\Delta$ be a Fuchsian group of signature $(0; 2,2,2,4)$ with canonical presentation $$\Delta=\langle \gamma_1, \gamma_2, \gamma_3, \gamma_4 : \gamma_1^2=\gamma_2^2=\gamma_3^2=\gamma_4^4=\gamma_1\gamma_2\gamma_3\gamma_4=1\rangle$$and, for each $n \geqslant 5, $ consider the group $G \cong (C_2 \times \mathbf{D}_{2^{n-3}}) \rtimes C_2$ of order $2^n$ with presentation in terms of generators $a,b,r,s$ and relations$$ r^{2^{n-3}}=s^2=(sr)^2=a^2=b^2=[s,b]=[r,b]=1, ara=r^{-1}, asa=sr, aba=br^{2^{n-4}}.$$

Note that the Riemann-Hurwitz formula is satisfied for a branched $2^n$-fold regular covering map from a compact Riemann surface of genus $1+2^{n-3}$ onto the projective line, ramified over three values marked with 2 and one value marked with 4. Thus, by virtue of Riemann's existence theorem, the existence of the desired family follows  after noticing that the correspondence $$\Delta \to G \,\, \mbox{ defined by }\,\, (\gamma_1, \gamma_2, \gamma_3, \gamma_4) \mapsto (s, bs, a, ab)$$is a surface-kernel epimorphism. Henceforth, we denote this family by $\mathcal{F}.$

\s

In order to determine an upper bound for the number of equisymmetric strata of $\mathcal{F}$, we have to determine an upper bound for the number of pairwise non-equivalent surface-kernel epimorphisms $\theta : \Delta \to G.$ For each such epimorphism $\theta$, we write$$g_i:=\theta(\gamma_i) \,\,\, \mbox{ for each } \,\,\, i =1,2,3,4,$$and, for the sake of simplicity, we identify $\theta$ with the $4$-uple  $\theta=(g_1, g_2, g_3, g_4).$
\s

We notice that: 
\begin{enumerate}
\item the elements of order four of $G$ are $abr^{l}$ and $w:=r^{2^{n-5}},$ and 
\item the involutions  of $G$ are $b, ar^l, z:=r^{2^{n-4}},  bz, sr^l$ and $bsr^l,$   
\end{enumerate}where $1 \leqslant l \leqslant 2^{n-3}.$

\s

{\bf Claim.} The central element $z$ is different from $g_1, g_2$ and $g_3.$

\s

Clearly, not three or two among $g_1, g_2, g_3$ can be equal to  $z.$ In addition, if one of them equals $z$, say $g_1=z,$ then, as $g_1g_2g_3$ must have order four, either $$g_2, g_3 \in \{ sr^{l_1}, bsr^{l_2} : 1 \leqslant l_j \leqslant 2^{n-3} \} \,\, \mbox{ or }\,\, g_2, g_3 \notin \{ sr^{l_1}, bsr^{l_2} : 1 \leqslant l_j \leqslant 2^{n-3}\}.$$

In the former case $g_4$ does not have order four, whilst in the latter one $s$ does not belong to the image of $\theta,$ contradicting its surjectivity.
\s

Similarly as argued before, we can see that the number of $g_i's$ that are of the form $sr^l$ or $bsr^l$ is exactly two. In addition, if $g_4$ equals $w$ then $\langle g_1, g_2, g_4 \rangle$ is a proper subgroup of $G$, showing that $\theta$ is not surjective. Thus, $\theta$ is of one of the following forms: $$(sr^{l_1}, sr^{l_2}, g_3, abr^{{l_3}}), \,\, (sr^{l_1}, bsr^{l_2}, g_3, abr^{{l_3}}) \,\, \mbox{or} \,\, (bsr^{l_1}, bsr^{l_2}, g_3, abr^{{l_3}})$$for some $1 \leqslant l_1, l_2, l_3 \leqslant 2^{n-3}.$ The fact that $g_1g_2g_3g_4=1$ implies that necessarily $\theta$ is $$ (sr^{l_1}, bsr^{l_2}, ar^{l_2+l_3-l_1}, abr^{{l_3}}) \,\, \mbox{ for some }\,\, 1 \leqslant l_1, l_2, l_3 \leqslant 2^{n-3}.$$

Note that, after applying an appropriate conjugation, we can assume $l_1=0$ or $l_1=1.$ Furthermore, by considering the action of the automorphism of $G$ given by $$r \mapsto r^{-1}, \,\, s \mapsto sr, \,\, a \mapsto a, \,\, b \mapsto b$$we obtain that $\theta$ is equivalent to $$\theta_{u,v}:=(s, bsr^{u}, ar^{v}, abr^{v-u}) \,\, \mbox{ where }\,\,1 \leqslant u,v, \leqslant 2^{n-3}.$$Thereby, the number of topologically non-equivalent  actions of $G$ on $C$ is as most $2^{2n-6}.$

\s

Following  \cite[Theorem 1]{singerman2}, the signature $(0; 2,2,2,4)$ is maximal; thus,  if $C$ lies in the interior of the family then its automorphism group agrees with $G.$ It is easy to verify that $G$ has exactly five conjugacy classes of subgroups of order two, and that among them only $K=\langle z \rangle$ is a normal subgroup. Consider the associated two-fold regular covering map given by the action of $K$ $$ \pi : C \to C_K,$$and notice that, independently of the equisymmetric stratum to which $C$ belongs (or, in other words, independently of the surface-kernel epimorphism $\theta_{u,v}$ representing the corresponding action), the covering $\pi$ ramifies over exactly $2^{n-2}$ values marked with 2. Thus, the Riemann-Hurwitz formula implies that $C_K$ is an elliptic curve and therefore $C$ is non-hyperelliptic. 

\s

If a compact Riemann surface $X$ belongs to $\mathcal{F}$ but does not belong to its interior then $G$ is strictly contained in the full automorphism group of $X$ (this is a general result that can be found, for instance, in \cite{Brou}).  Now, as the complex dimension of the family $\mathcal{F}$ is one, it follows that the signature of the action of $\mbox{Aut}(X)$ on $X$ must be triangle; namely, of the form $(0; t_1, t_2, t_3).$ Note that there are finitely many possibilities for $t_1, t_2, t_3$ and, in turn, to each of these possible signatures correspond at most finitely many Riemann surfaces. Thus, the family contains at most finitely many surfaces that do not belong to its interior.

\s

We now proceed to prove the announced isogeny decomposition of $JC$ for each $C$ in the family $\mathcal{F}$. Let us consider the normal subgroup $N$ of $G$ given by $$\langle r,s,b: r^{2^{n-3}}=s^2=(sr)^2=b^2=1, [b,s]=[b,r]=1 \rangle$$ and the complex irreducible representation of $N$ given by $$r \mapsto \left( \begin{smallmatrix}
\omega & 0 \\
0 & \bar{\omega}
\end{smallmatrix} \right), \,\,\,\, s \mapsto \left( \begin{smallmatrix}
0 & 1 \\
1 & 0
\end{smallmatrix} \right), \,\,\,\, b \mapsto \left( \begin{smallmatrix}
1 & 0 \\
0 & 1
\end{smallmatrix} \right) \,\, \mbox{ where }  \omega \mbox{ is a }2^{n-3}\mbox{-th primitive root of unity.}$$

This representation induces the complex representation $V$ of $G$ given by $$r \mapsto \left( \begin{smallmatrix}
\omega & 0 & 0 & 0\\
0 & \bar{\omega} & 0 & 0\\
0 & 0 & \bar{\omega} & 0\\
0 & 0 &0 & \omega
\end{smallmatrix} \right), \,\, s \mapsto \left( \begin{smallmatrix}
0 & 1 & 0 & 0\\
1 & 0 & 0 & 0\\
0 & 0 & 0 & \omega\\
0 & 0 & \bar{\omega} & 0
\end{smallmatrix} \right), \,\, b \mapsto  \left( \begin{smallmatrix}
1 & 0 & 0 & 0\\
0 & 1 & 0 & 0\\
0 & 0 & -1 & 0\\
0 & 0 &0 & -1
\end{smallmatrix} \right), \,\, a  \mapsto  \left( \begin{smallmatrix}
0 & 0 & 1 & 0\\
0 & 0 & 0 & 1\\
1 & 0 & 0 & 0\\
0 & 1 &0 & 0
\end{smallmatrix} \right)$$which is,  by  \cite[Theorem 6.11]{issac},  irreducible. In addition, as $V$ is constructed from a  complex irreducible representation of a dihedral group, it is easy to infer that its Schur index is 1. Note  the character field  of $V$ is $\mathbb{Q}(\omega + \bar{\omega});$ this is an extension of $\mathbb{Q}$ of degree $$\tfrac{1}{2}\varphi(2^{n-3})=2^{n-5}$$where $\varphi$ is the Euler function. We denote by $W_2$ the rational irreducible representation of $G$ associated to $V$ and by $W_1$ the rational irreducible representation $G$ given by $$r \mapsto 1, \,\,\, s \mapsto -1, \,\,\, b \mapsto 1, \,\,\, a \mapsto -1.$$

Then, as explained in \S \ref{jacos}, there is an abelian subvariety $P$ of $JC$ such that \begin{equation} \label{deco}JC \sim B_{W_1} \times B_{W_2}^4 \times P,\end{equation}where $B_{W_j}$ is the factor associated to $W_j$ in the group algebra decomposition of $JC$ with respect to $G.$ 
As the action of $G$ on $C$ is determined by $\theta_{u,v}$ for some $u,v \in \{1, \ldots, 2^{n-3}\},$ we can apply the equation \eqref{uuaa}  to notice that, independently of the choice of $u$ and $v,$ the following equalities hold:   $$\dim(B_{W_1})=1 \,\, \mbox{ and } \,\, \dim(B_{W_2})=2^{n-5}.$$ Then, by considering dimensions is the relation \eqref{deco}, one sees that $$\dim (JC) =1+2^{n-3}=1+ 4( 2^{n-5} )+ \dim P \,\, \mbox{ and therefore }\,\, P=0.$$Now, we consider the induced isogeny \eqref{indjaco}  (with $H=\langle b \rangle$ and $H=\langle s \rangle$) to obtain  $$JC_{\langle b \rangle} \sim B_{W_1} \times B_{W_2}^2\,\, \mbox{ and } \,\, JC_{\langle s \rangle}  \sim B_{W_2}^{2}$$

The previous two isogenies together with isogeny \eqref{deco} permits us to conclude that $$JC \sim JC_{\langle b \rangle} \times JC_{\langle s \rangle}$$as claimed. Finally, is a similar way, we consider the induced isogeny \eqref{indjaco}  with $H=\langle r \rangle$ to obtain  that $JC_{\langle r \rangle}$ and $B_{W_1}$ are isogenous and, consequently, $JC$ contains an elliptic curve isogenous to $JC_{\langle r \rangle}$.

\section{Proof of Theorem \ref{ppp}} \label{s5}

Let $d \geqslant 1, g \geqslant 2$ be integers and let $p \geqslant 3$ be a prime number. Let $G$ be a $p$-group of automorphisms of a complex $d$-dimensional family of compact Riemann surfaces $C$ of genus $g$ and assume the signature of the action of $G$ on $C$ to be $\sigma=(h; m_1, \ldots, m_l).$

The hyperbolic area $\mu$ of a fundamental region of a Fuchsian group of signature $\sigma$ satisfies  \begin{displaymath} \mu \geqslant 2\pi[d+1-\tfrac{d+3}{p}+h(\tfrac{3}{p}-1)] \geqslant \left\{ \begin{array}{ll}
\tfrac{4}{3}d\pi & \textrm{if $p=3$}\\
2\pi[d+1-\tfrac{d+3}{p}+\hat{h}(\tfrac{3}{p}-1)] & \textrm{if $p \geqslant 5$}
  \end{array} \right.
\end{displaymath}where $\hat{h}$ is the largest possible genus of the quotient $C_G.$ Note that $\hat{h}=[\tfrac{d}{3}+1].$

\s

Assume $p=3.$ The Riemann-Hurwitz formula ensures that $$2(g-1) =|G| \tfrac{\mu}{2 \pi} \geqslant  \tfrac{2}{3}d|G| \,\, \iff \,\, |G| \leqslant M_{3,d}$$as claimed in (1). Now, if  we suppose that the order of $G$ equals $M_{3,d}$ then, by the Riemann-Hurwitz formula, we easily obtain that $$\Sigma_{i=1}^l \tfrac{1}{m_i}=\tfrac{l}{3} \,\, \mbox{ and, consequently, each } m_i=3.$$ Note that there is no restriction on $l.$ Thus, $\sigma=\sigma_{3,d,h}$ for some $h \in \{0, \ldots, \hat{h}\}.$ The {\it only if} part of (2) is a direct computation.

\s

Assume $p \geqslant 5.$ Then\begin{displaymath} \mu \geqslant  \left\{ \begin{array}{ll}
\tfrac{4}{3} \pi d& \textrm{if $d \equiv 0 \mbox{ mod } 3$}\\
\tfrac{4}{3} \pi d+2 \pi(\tfrac{1}{3}-\tfrac{1}{p})& \textrm{if $d \equiv 1 \mbox{ mod } 3$}\\
\tfrac{4}{3} \pi d+4 \pi(\tfrac{1}{3}-\tfrac{1}{p}) & \textrm{if $d \equiv 2 \mbox{ mod } 3.$}
  \end{array} \right.
\end{displaymath} In other words, if $\lambda_d$ is the smallest non-negative representative of $d$ modulo 3 then $$2(g-1)=|G|\tfrac{\mu}{2 \pi} \geqslant |G|(\tfrac{2}{3}d+\lambda_d(\tfrac{1}{3}-\tfrac{1}{p})) \iff |G| \leqslant M_{p,d},$$as claimed in (3). If we now assume that the order of $G$ equals $M_{p,d}$ then  \begin{equation} \label{tele}
\Sigma_{i=1}^l
 \tfrac{1}{m_i}=\tfrac{l}{3}-\lambda_d(\tfrac{1}{3}-\tfrac{1}{p}).\end{equation}
 
\begin{enumerate}
\item If $d \equiv 0 \mbox{ mod } 3$ then \eqref{tele} turns into $\Sigma_{i=1}^l
 \tfrac{1}{m_i}=\tfrac{l}{3}$ and $l=0.$  Thus, $$\sigma=(\tfrac{d+3}{3}; -)=\sigma_{p,d}$$
\item If $d \equiv 1 \mbox{ mod } 3$ then \eqref{tele} turns into $\Sigma_{i=1}^l
 \tfrac{1}{m_i}=\tfrac{l}{3}-\tfrac{1}{3}+\tfrac{1}{p}$ and $l=1.$ Thus, $$\sigma=(\tfrac{d+2}{3}; p)=\sigma_{p,d}$$
\item If $d \equiv 2 \mbox{ mod } 3$ then \eqref{tele} turns into $\Sigma_{i=1}^l
 \tfrac{1}{m_i}=\tfrac{l}{3}-\tfrac{2}{3}+\tfrac{2}{p}$ and $l=2.$ Thus, $$\sigma=(\tfrac{d+1}{3}; p,p)=\sigma_{p,d}$$
\end{enumerate}

The {\it only if} part of (4) is a direct computation.

\section{Proof of Theorem \ref{pexistencia}} \label{s6}
Let $p \geqslant 3$ and let $\Delta$ be a Fuchsian group of signature $(1;p)$ with canonical presentation $$\Delta=\langle \alpha_1, \beta_1, \gamma_1: \alpha_1\beta_1\alpha_1^{-1}\beta_1^{-1}\gamma_1=\gamma_1^p=1   \rangle$$and, for each $n \geqslant 5, $ consider the group  $G \cong C_{p^{n-1}} \rtimes_p C_p$ of order $p^n$ with presentation $$\langle a,b : a^{p^{n-1}}=b^{p}=1, bab^{-1}=a^{r} \rangle,$$where $r=p^{n-2}+1.$ Observe that $r^p \equiv 1 \mbox{ mod } p^{n-1}$ and $r^k \not\equiv 1 \mbox{ mod } p^{n-1}$ for $1 \leqslant k \leqslant  p-1.$

Note that the Riemann-Hurwitz formula is satisfied for a branched $p^n$-fold regular covering map from a compact Riemann surface of genus $1+\tfrac{(p-1)p^{n-1}}{2} $ onto a Riemann surface of genus one, ramified over one value marked with $p$. Thus, by virtue Riemann's existence theorem, the existence of the family follows  after noticing that the rule $$\Delta \to G \,\, \mbox{ defined by }\,\, (\alpha_1, \beta_1, \gamma_1) \mapsto (a,b,a^{p^{n-2}})$$is a surface-kernel epimorphism.  Henceforth, we denote this family by $\mathcal{G}.$

\s

We now proceed to prove that there are exactly $p-1$  pairwise non-equivalent surface-kernel epimorphisms $\theta: \Delta\to G.$ For each such epimorphism $\theta$, we write$$x:=\theta(\alpha_1), \,\, y=\theta(\beta_1) \,\, \mbox{ and }\,\, z=\theta(\gamma_1)$$and, for the sake of simplicity, we identify $\theta$ with the $3$-uple  $\theta=(x,y,z).$ Note that $$x=a^l b^k \,\, \mbox{ and }\,\, y=a^sb^m \,\, \mbox{ for some }\,\, 1 \leqslant l,s \leqslant p^{n-1} \,\, \mbox{ and } \,\, 1 \leqslant k,m \leqslant p. $$

If $k \neq 0$ and $u=-mk',$ where $k'$ is the inverse of $k$ in the field of $p$ elements, then the transformation $A_{1,u}$ (see \S \ref{strati}) shows that we can assume, up to equivalence, that  \begin{equation} \label{virus}x=a^l b^k \,\, \mbox{ and } \,\,  y=a^s.\end{equation} 
On the other hand, if $k=0$ then \begin{equation} \label{virus2}x=a^l  \,\, \mbox{ and } \,\,  y=a^sb^m,\end{equation}and the transformation $A_{2,-1}\circ A_{1,1}$ 
shows that  \eqref{virus} and \eqref{virus2} are equivalent. Now, in \eqref{virus2} one sees that if $l$ and $p^{n-1}$ are not coprime then $\theta$ is not surjective. Thus, after sending $a$ to an appropriate power of it, we can be assume $l=1.$ Then  \begin{equation} \label{virus3} x=a \,\, \mbox{ and }\,\, b=a^sb^m \,\, \mbox{ where } \,\, m \neq 0.\end{equation}Now, if we set $v=-sr^{-m}$ then we apply 
$A_{1,v}$ to \eqref{virus3} to ensure that $\theta$ is equivalent to $$\theta_m=(a, b^m, a^{r^m-1}) \,\, \mbox{ for some }\,\, 1 \leqslant m \leqslant p-1.$$

The result follows after noticing that $\theta_m$ and $\theta_{m'}$ are non-equivalent if $m \neq m'.$

\s

\s

Note that $K=\langle a^{p^{n-2}} \rangle$ is a cyclic group of order $p$ and that, independently of the equisymmetric stratum to which $C$ belongs, the associated regular covering map $$C \to C_K$$ramifies over $p^{n-1}$ values marked with $p.$ It follows that the quotient Riemann surface $C_K$ has genus one; thus, $C$ is an elliptic-$p$-gonal Riemann surface. Due to the explicitness of the family, one can easily see that $K$ is the unique group of automorphisms of $C$ providing the elliptic-$p$-gonal structure (c.f. \cite[Theorem 1.3]{Sch4} and also \cite{GroHi} and \cite{GroWW}).

\s

According to \cite[Theorem 1]{singerman2},  the action of $G$ on each $C$ in $\mathcal{G}$ might be extended to only an action of a group of order $2p^n$ acting on $C$ with signature $\sigma'=(0; 2,2,2,2p).$ 

\s

{\bf Claim.} Such extension is not possible in our case. 

\s

To prove the claim we shall proceed by contradiction; namely, we assume that:

\begin{enumerate}
\item there is a group $G'$ of order $2q^n$ with a subgroup isomorphic to $G,$ and that 
\item there is a surface-kernel epimorphism $\Delta' \to G',$ where $\Delta'$ is a Fuchsian group of signature $\sigma'.$ 
\end{enumerate}

By the classical Schur-Zassenhaus theorem, we can ensure that $$G' \cong G \rtimes C_2 \,\, \mbox{ with } \,\, C_2=\langle t : t^2 =1 \rangle.$$Observe that $C_2$ must act on $G$ with order 2, because of the direct product $G \times C_2$ cannot be generated by three involutions. Thus, by considering an automorphism of $G$ that sends $a$ to an appropriate power of it and after some routine computations, one can  see that, up to an  isomorphism of $G$, the action of $C_2$ on $G$ is given by$$tat=a^{-1} \,\, \mbox{ and } \,\, tbt = b.$$In particular, the involutions of $G'$ are of the form $ta^k$ for $0 \leqslant k \leqslant p-1.$ However,  three of them cannot generate $G',$ contradicting the surjectivity of $\theta.$ This proves the claim.

\s

As observed in the proof of Theorem \ref{nil2}, the surface $C$ belongs to the interior of the family $\mathcal{G}$ (and therefore for all up to possibly finitely many exceptions) if and only if $G$ is the full automorphism group of it (see, for instance, \cite{Brou}).
  
\s

We now proceed to decompose the Jacobian variety $JC$ of each $C$ in the family $\mathcal{G}.$ 

We apply the method of {\it little groups} of Wigner--Mackey (see, for example, \cite[p. 62]{Serre}), to guarantee the irreducibility of the complex representation $V$ of $G$  given by $$a \mapsto \mbox{diag}( \omega, \omega^r, \ldots, \omega^{r^{p-1}}) \,\, \mbox{ and } \,\, b \mapsto\left( \begin{smallmatrix}
0 & 1 & 0 & \ldots & 0\\
0 & 0 & 1 & \cdots  & 0\\
\, & \, & \,  & \ddots & \, \\
0 & 0 & 0 & \ldots & 1\\
1 & 0 & 0 & \cdots & 0
\end{smallmatrix} \right)$$where $\omega$ is a $p^{n-1}$-th primitive root of unity. We notice that the character field of $V$ is $\mathbb{Q}(\omega^p)$, which is an extension of $\mathbb{Q}$ of degree $$\varphi(p^{n-2})=p^{n-3}(p-1),$$where $\varphi$ is the Euler function.  We recall that $p$-groups with $p \geqslant 3$ only possess representations with Schur index 1 (see, for example,   \cite[Theorem 41.9]{Reiner}).

We denote by $W$ the rational irreducible representation of $G$ associated to $V$. Then, as explained in \S \ref{jacos}, there is an abelian subvariety $Q$ of $JC$ such that \begin{equation} \label{deco2}JC \sim E \times A^p \times Q,\end{equation}where $E$ is an elliptic curve isogenous to $JC_G$ and $A$ is the factor associated to $W$ in the group algebra decomposition of $JC$ with respect to $G.$ Now, as the action of $G$ on $C$ is determined by $\theta_{m}$ for some $1 \leqslant m \leqslant p-1,$ we apply the equation \eqref{uuaa} to notice that, independently of the choice of $m,$  the following equality holds: $$ \dim(A)=p^{n-3}(p-1) \cdot \tfrac{1}{2}(p-0)=\tfrac{(p-1)p^{n-2}}{2}.$$ Finally, by considering dimensions in the relation \eqref{deco2}, one concludes that $Q=0$ and the desired decomposition of $JC$ is obtained.

\s

\subsection*{Acknowledgements} The author is very grateful to the referee for his/her  valuable comments and suggestions.


\begin{thebibliography}{9}


%
%
%
\bibitem{bci} 
{\sc G. Bartolini, A. F. Costa, and M. Izquierdo,} 
{\em On the orbifold structure of the moduli space of Riemann surfaces of genera four and five.}
Rev. R. Acad. Cienc. Exactas Fis. Nat. Ser. A Mat. RACSAM {\bf 108} (2014), no. 2, 769--793. 





\bibitem{BJ} 
{\sc M. V. Belolipetsky and G. A. Jones,} 
{\em Automorphism groups of Riemann surfaces of genus $p + 1$, where $p$ is prime.}
Glasg. Math. J. {\bf 47} (2005), no. 2, 379--393.





\bibitem{Brou}
{\sc{S. A. Broughton,}}  {\em{Classifying finite groups actions on surfaces of low genus,}} J. Pure Appl. Algebra {\bf 69} (1990), no. 3, 233--270.

%

\bibitem{b}
{\sc S. A. Broughton},
 { \em The equisymmetric stratification of the moduli space and the Krull dimension of mapping class groups,} Topology Appl. {\bf 37} (1990), no. 2, 101--113.
 
\bibitem{BCI}
{\sc{E. Bujalance, A. F. Costa and M. Izquierdo,}}  {\em{On Riemann surfaces of genus g with 4g automorphisms,}} Topology Appl. {\bf 218} (2017) 1--18.

\bibitem{CLR}
{\sc A. Carocca, H. Lange and R. E. Rodriguez},
 { \em \'Etale double covers of cyclic p-gonal covers,}  J. Algebra {\bf 538} (2019), 110--126.


\bibitem{d1}
{\sc A. Carocca, S. Recillas and R. E. Rodr\'iguez},
 { \em Dihedral groups acting on Jacobians,} Contemp. Math. {\bf 311} (2011), 41--77.

\bibitem{CRC}
{\sc A. Carocca and S. Reyes-Carocca},
 { \em Riemann surfaces of genus $1+q^2$ with $3q^2$ automorphisms,} Preprint, arXiv:1911.04310v1.
 
\bibitem{cr}
{\sc A. Carocca and R. E. Rodr\'iguez,}
{\em Jacobians with group actions and rational idempotents.}
J. Algebra \textbf{306} (2006), no. 2, 322--343.

\bibitem{Che}
{\sc B. P. Chetiya}, {\em On genuses of compact Riemann surfaces admitting solvable automorphism groups}, Indian J. Pure Appl. Math. {\bf 12} (1981), 1312--1318.


\bibitem{ChP}
{\sc B. P. Chetiya and K. Patra}, {\em  On metabelian groups of automorphisms of compact Riemann surfaces,} J. London Math. Soc. {\bf 33} (1986), 467--472.


\bibitem{CI}
{\sc A. F. Costa and M. Izquierdo,}
{\em One-dimensional families of Riemann surfaces of genus $g$ with $4g + 4$ automorphisms,} Rev. R. Acad. Cienc. Exactas F\'is. Nat. Ser. A Mat. RACSAM {\bf 112} (2018), no. 3, 623--631.

\bibitem{Do}
{\sc R. Donagi and E. Markman}, {\em Spectral covers, algebraically completely integrable, Hamiltonian systems, and moduli of bundles}, in: Integrable Systems and Quantum Groups, Montecatini Terme, 1993, in: Lecture Notes in Math., vol. 1620, Springer, Berlin, 1996, pp. 1--119.

\bibitem{Gro}
{\sc G. Gromadzki,} {\em Maximal groups of automorphisms of compact Riemann surfaces in various classes of finite groups}, Rev. Real Acad. Cienc. Exact. Fis. Natur.  {\bf 82} no. 2 (1988), 267--275.

\bibitem{Gro4}
{\sc G. Gromadzki,} {\em Metabelian groups acting on compact Riemann surfaces}, Rev. Mat. Univ. Complut. Madrid {\bf 8} no. 2 (1995), 293--305.

\bibitem{Gro2}
{\sc G. Gromadzki,} {\em On soluble groups of automorphism of Riemann surfaces}, Canad. Math. Bull. {\bf 34} (1991), 67--73.


\bibitem{Gro3}
{\sc G. Gromadzki and C. Maclachlan,} {\em Supersoluble groups of automorphisms of compact Riemann surfaces}, Glasgow Math. J. {\bf 31} (1989), 321--327.

\bibitem{GroHi}
{\sc G. Gromadzki and R. A. Hidalgo}, {\em  On prime Galois coverings of tori,} preprint.

\bibitem{GroWW}
{\sc G. Gromadzki, A. Weaver and A. Wootton}, {\em  On gonality of Riemann surfaces,} Geom. Dedicata {\bf 149} (2010), 1--14. 

\bibitem{Harvey1}
{\sc J. Harvey,}
{\em Cyclic groups of automorphisms of a compact Riemann surface}. 
Q. J. Math. {\bf 17}, (1966), 86--97.

\bibitem{Harvey}
{\sc J. Harvey,}
{\em On branch loci in Teichm\"{u}ller space},
Trans. Amer. Math. Soc. {\bf 153} (1971), 387--399.

\bibitem{nos}
{\sc R. A. Hidalgo, L. Jim\'enez, S. Quispe and S. Reyes-Carocca,} {\em Quasiplatonic curves with symmetry group $\mathbb{Z}_2^2 \rtimes \mathbb{Z}_m$ are definable over $\mathbb{Q}$,} Bull. London Math. Soc. {\bf 49} (2017) 165--183.

\bibitem{issac}
{\sc I. M. Isaacs,}
{\em Character theory of finite groups}, Academic Press, (1976).

\bibitem{IJR}
{\sc M. Izquierdo, L. Jim\'enez, A. Rojas,}
{\em Decomposition of Jacobian varieties of curves with dihedral actions via equisymmetric stratification},   Rev. Mat. Iberoam. {\bf 35}, No. 4 (2019), 1259--1279.

\bibitem{IJRC}
{\sc M. Izquierdo, G. A. Jones and S. Reyes-Carocca,}
{\em Groups of automorphisms of Riemann surfaces and maps of genus $p+1$ where $p$ is prime}, To appear in Ann. Acad. Sci. Fenn. Math. (2021), arXiv:2003.05017.

\bibitem{IRC}
{\sc M. Izquierdo and S. Reyes-Carocca,}
{\em A note on large automorphism groups of compact Riemann surfaces,}
 J. Algebra {\bf 547} (2020), 1--21.

\bibitem{IRCR}
{\sc M. Izquierdo, S. Reyes-Carocca and A. M. Rojas,}
{\em On families of Riemann surfaces with automorphisms,} J. Pure Appl. Algebra {\bf 224} no. 10, 106704 (2021).

\bibitem{KR}
{\sc E. Kani and M. Rosen,} {\em Idempotent relations and factors of Jacobians},
Math. Ann. {\bf 284} (1989), 307--327.

\bibitem{K1}
{\sc R. S. Kulkarni,}
{ \em A note on Wiman and Accola-Maclachlan surfaces}.
Ann. Acad. Sci. Fenn., Ser. A 1 Math. {\bf 16} (1) (1991)
83--94.

\bibitem{l-r}
{\sc H. Lange and S. Recillas,}
{ \em Abelian varieties with group actions}.
J. Reine Angew. Mathematik, \textbf{575} (2004) 135--155.

\bibitem{McB2}
{\sc A. MacBeath,}
{\em Action of automorphisms of a compact Riemann surface on the first homology group}, Bull. London Math. Soc. {\bf 5} (1973), 103--108.

\bibitem{MMBB}
{\sc A. MacBeath,}
{\em Residual Nilpotency of Fuchsian groups,} Illinois Journal of Mathematics {\bf 28} no. 2 (1984), 299--311.

\bibitem{Mac2}
{\sc C. Maclachlan,} {\em Abelian groups of automorphisms of compact Riemann surfaces}, Proc. London Math. Soc. {\bf 15} (1965), 699--712.

\bibitem{PA}
{\sc J. Paulhus and A. M. Rojas,}
{ \em Completely decomposable Jacobian varieties in new genera},  Experimental Mathematics {\bf 26} (2017), no. 4, 430--445.

\bibitem{d3}
{\sc S. Recillas and R. E. Rodr\'iguez,} {\em Jacobians and representations of $S_3$}, Aportaciones Mat. Investig. {\bf 13}, Soc. Mat. Mexicana, M\'exico, 1998.

\bibitem{Reiner}
{\sc I. Reiner,} {\em Maximal orders}, London Math. Soc. Monogr. (N. S.), vol. {\bf 28}, Oxford Univ. Press, Oxford 2003.

\bibitem{RR1}
{\sc S. Reyes-Carocca,} {\em On Riemann surfaces of genus $g$ with $4g-4$ automorphisms}, Israel J. Math. {\bf 237}, 415--436 (2020).

\bibitem{yo-racsam}
{\sc S. Reyes-Carocca,} {\em On pq-fold regular covers of the projective line,}
Rev. R. Acad. Cienc. Exactas Fís. Nat. Ser. A Mat. {\bf 115}, 23 (2021).

\bibitem{yojpaa}
{\sc S. Reyes-Carocca,} {\em On the one-dimensional family of Riemann surfaces of genus $q$ with $4q$ automorphisms}, J. of Pure Appl. Algebra 223, no. {\bf 5} (2019), 2123--2144.




\bibitem{kanirubiyo}
{\sc S. Reyes-Carocca and R. E. Rodr\'iguez,} {\em A generalisation of Kani-Rosen decomposition theorem for Jacobian varieties}, Ann. Sc. Norm. Super. Pisa Cl. Sci. (5) {\bf 19} (2019), no. 2, 705--722.

\bibitem{RCR}
{\sc S. Reyes-Carocca and R. E. Rodr\'iguez,} {\em On Jacobians with group action and coverings}, Math. Z. (2020) {\bf 294}, 209--227.

\bibitem{Ri}
{\sc J. Ries}, {\em The Prym variety for a cyclic unramified cover of a hyperelliptic curve}, J. Reine Angew. Math. {\bf 340} (1983) 59--69.

\bibitem{yoibero}
{\sc A. M.   Rojas}, {\em Group actions on Jacobian varieties}, Rev. Mat. Iber. {\bf 23} (2007), 397--420.

\bibitem{Sch2}
{\sc A. Schweizer},  {\em Metacyclic groups as automorphism groups of compact Riemann
surfaces}, Geom. Dedicata {\bf 190} (2017), 185--197.

\bibitem{Sch4}
{\sc A. Schweizer},  {\em  On the uniqueness of $(p,h)$-gonal automorphisms of Riemann surfaces.} Arch. Math. (Basel) {\bf 98} (2012), no. 6, 591--598. 

\bibitem{Sch1}
{\sc A. Schweizer},  {\em Several types of solvable groups as automorphism groups of compact Riemann surfaces}, 	arXiv:1701.00325.

\bibitem{Serre}
{\sc J. P. Serre}, {\em Linear Representations of Finite Groups}, Graduate Texts in Mathematics. {\bf 42} Springer-Verlag, New York.

\bibitem{singerman2}
{\sc D. Singerman}, {\em Finitely maximal Fuchsian groups}, J. London Math. Soc. (2)  {\bf 6}, (1972), 29--38.

\bibitem{singerman}
{\sc D. Singerman}, {\em Subgroups of Fuchsian groups and finite permutation groups}, Bull. London Math. Soc.  {\bf 2}, (1970), 319--323.

\bibitem{Z2}
{\sc R. Zomorrodian}, {\em  Bounds for the order of supersoluble automorphism groups of Riemann surfaces}, Proc. Amer. Math. Soc. {\bf 108} no. 3 (1990), 587--600.

\bibitem{Z3}
{\sc R. Zomorrodian}, {\em  Classification of p-groups of automorphisms of Riemann surfaces and their lower central series,} Glasgow Math. J. {\bf 29} (1987), 237--244.

\bibitem{Z1}
{\sc R. Zomorrodian}, {\em Nilpotent automorphism groups of Riemann surfaces,}   Trans. Amer. Math. Soc. {\bf 288} (1985), no. 1, 241--255.

\bibitem {Wi}
{\sc  A. Wiman}, 
{\em  \"{U}ber die hyperelliptischen Curven und diejenigen von Geschlechte p - Jwelche eindeutige Tiansformationen in sich zulassen.} \textit{Bihang till K. Svenska Vet.-Akad. Handlingar, Stockholm} \textbf{21}
(1895-6) 1-28.



\end{thebibliography}
\end{document}